\documentclass[12pt,reqno]{amsart}

\usepackage{amssymb}

\def\fg{{\mathfrak g}}

\def\fm{{\mathfrak m}}

\def\BN{{\mathbb N}}

\def\BN{{\mathbb N}}

\def\b{\operatorname{b}}
\def\fg{\operatorname{fg}}

\def\opp{\operatorname{opp}}

\def\CMreg{\operatorname{CMreg}}

\def\D{\mbox{\sf D}}

\def\Ext{\operatorname{Ext}}
\def\Extreg{\operatorname{Ext.\!reg}}

\def\Gr{\mbox{\sf Gr}}

\def\h{\operatorname{h}}
\def\HGammam{\operatorname{H}_{\fm}}
\def\HGammamopp{\operatorname{H}_{\fm^{\opp}}}
\def\Hom{\operatorname{Hom}}

\def\LTensor{\stackrel{\operatorname{L}}{\otimes}}

\def\RGammam{\operatorname{R}\!\Gamma_{\fm}}
\def\RGammamopp{\operatorname{R}\!\Gamma_{{\fm}^{\opp}}}

\def\RHom{\operatorname{RHom}}

\def\Tor{\operatorname{Tor}}

\def\Tot{\operatorname{Tot}}

\numberwithin{equation}{part}

\newtheorem{Lemma}{Lemma}[section]
\newtheorem{Theorem}[Lemma]{Theorem}
\newtheorem{Proposition}[Lemma]{Proposition}
\newtheorem{Corollary}[Lemma]{Corollary}
\theoremstyle{definition}
\newtheorem{Definition}[Lemma]{Definition}
\newtheorem{Setup}[Lemma]{Setup}

\newtheorem{Observation}[Lemma]{Observation}

\begin{document}

\title[Linear free resolutions]
{Linear free resolutions over non-commutative algebras}

\author{Peter J\o rgensen}
\address{Danish National Library of Science and Medicine, N\o rre
All\'e 49, 2200 K\o \-ben\-havn N, DK--Denmark}
\email{pej@dnlb.dk, www.geocities.com/popjoerg}
%\thanks{Date: \today. A thank you would go here}

%\author{Next author goes here}
%\address{Next author's postal address goes here}
%\email{Next author's mail address goes here}

%%  2000 AMS subject classification:
%%  16E05: Syzygies, resolutions, complexes
%%  16E30: Homological functors on modules (Tor, Ext, etc.)
%%  16W50: Graded rings and modules

\begin{abstract} 
The main result of this paper is that over a non-commutative Koszul
algebra, high truncations of finitely generated graded modules have
linear free resolutions.
\end{abstract}

\maketitle

\setcounter{section}{-1}
\section{Introduction}
\label{sec:introduction}

\noindent
Eisenbud and Goto, Avramov and Eisenbud, and this author have all
studied whether high truncations of finitely generated graded modules
over graded algebras have linear free resolutions, see
\cite{AE}, \cite{EG}, and \cite{PJCM}.

The original study of this took place over polynomial algebras.  The
main result is \cite[prop.\ p.\ 89]{EG}: If $M$ is a finitely
generated graded module over $k[X_1, \ldots, X_t]$ with $k$ a field,
then for large $s$, the minimal free resolution of the degree shifted
truncation $M_{\geq s}(s)$ is linear. That is, the $m$'th module in
the minimal free resolution has all its generators placed in degree
$m$.

This was later extended in \cite[cor.\ 2]{AE} to commutative Koszul
algebras, and in \cite[thm.\ 2.6]{PJCM} to non-commutative AS-regular
algebras, which are algebras with good homological behaviour
generalizing that of polynomial algebras.

This paper proves a common extension of \cite[cor.\ 2]{AE} and
\cite[thm.\ 2.6]{PJCM}: If $A$ is a non-commutative Koszul algebra
satisfying a few weak conditions given in setup \ref{set:blanket},
then for any finitely generated graded module $M$ and large $s$, the
minimal free resolution of $M_{\geq s}(s)$ is linear. This is theorem
\ref{thm:linear} below.

Along the way, I prove theorems \ref{thm:main1} and \ref{thm:main2}
and corollaries \ref{cor:regs_coincide} and \ref{cor:Extreg_finite}
which show that the two competing definitions of Castelnuovo-Mumford
regularity of graded modules given in \cite{AE}, respectively
\cite{EG} and \cite{PJCM}, are in fact closely related.

\begin{Setup}
\label{set:blanket}
Throughout the paper, $k$ is a field, and $A$ is a noetherian
$\BN$-graded connected $k$-algebra which has a balanced dualizing
complex. 
\end{Setup}

See \cite{PJGorHom} and \cite{PJCM} for generalities on the theory of
graded algebras, and \cite{VdBExist} and \cite{YDual} for information
on dualizing complexes.  My notation is mostly standard, but I do want
to give a few keywords:

The opposite algebra of $A$ is denoted $A^{\opp}$, and
$A$-right-modules are identified with $A^{\opp}$-left-modules.

The abelian category of graded $A$-left-modules and graded
homomorphisms of degree zero is denoted $\Gr\: A$.  The derived
category of $\Gr\: A$ is denoted $\D(\Gr\: A)$.  If $X$ is in
$\D(\Gr\: A)$, then $\h^m\!X$ denotes the $m$'th cohomology module of
$X$.  The derived category $\D(\Gr\: A)$ has full subcategories
$\D^-(\Gr\: A)$ consisting of complexes $X$ with $\h^m\!X = 0$ for $m$
large positive, $\D^+(\Gr\: A)$ consisting of complexes $X$ with
$\h^m\!X = 0$ for $m$ large negative, and $\D^{\b}_{\fg}(\Gr\: A)$
consisting of complexes $X$ with $\h^m\!X = 0$ for $m$ large positive
or negative and each $\h^m\!X$ finitely generated.

The derived functors of $\Hom_A$ and $\otimes_A$ are denoted $\RHom_A$
and $\LTensor_A$.  Section functors of graded $A$-left- and graded
$A$-right-modules are denoted $\Gamma_{\fm}$ and
$\Gamma_{\fm^{\opp}}$, and their derived functors are denoted
$\RGammam$ and $\RGammamopp$.  These give rise to local cohomology
functors by $\HGammam^m = \h^m \RGammam$ and $\HGammamopp^m = \h^m
\RGammamopp$.  The Matlis duality functor is denoted $(-)^{\prime}$
and defined on graded $A$-modules by $(M^{\prime})_p =
\Hom_k(M_{-p},k)$.  Matlis duality exchanges graded $A$-left- and
graded $A$-right-modules, and is exact and therefore well-defined on
derived categories.

\section{Background results}
\label{sec:background}

\begin{Proposition}
\label{prp:RHom_isomorphism}
For $X$ and $Y$ in $\D^{\b}_{\fg}(\Gr\: A)$ there is a natural
isomorphism 
\[
  \RHom_A(\RGammam X,Y) \cong \RHom_A(X,Y).
\]
\end{Proposition}

\begin{proof}
First observe that there are natural isomorphisms
\begin{align*}
  \RHom_A(\RGammamopp A,Y) 
  & \cong \RHom_A(\RGammamopp A,Y^{\prime \prime}) \\
  & \stackrel{\rm (a)}{\cong} (Y^{\prime} \LTensor_A \RGammamopp A)^{\prime} \\
  & \stackrel{\rm (b)}{\cong} \RGammamopp(Y^{\prime})^{\prime} \\
  & \stackrel{\rm (c)}{\cong} Y^{\prime \prime} \\
  & \cong Y.
\end{align*} 
Here (a) is by \cite[thm.\ 1.5]{PJGorHom} and (b) is by
\cite[thm.\ 1.6]{PJGorHom}, while (c) can be seen as follows:
Let $F$ be a free resolution of $Y$ consisting of finitely generated
free modules. It is then easy to see that $F^{\prime}$ is an injective
resolution of $Y^{\prime}$. As $F$ consists of finitely generated free
modules, $F^{\prime}$ consists of torsion graded injective modules, so
$\Gamma_{\fm^{\opp}}(F^{\prime}) \cong F^{\prime}$ whence
$\RGammamopp(Y^{\prime}) \cong \Gamma_{\fm^{\opp}}(F^{\prime}) \cong
F^{\prime} \cong Y^{\prime}$.

Now compute:
\begin{align*}
  \RHom_A(\RGammam X,Y)
  & \stackrel{\rm (d)}{\cong} \RHom_A(\RGammam(A) \LTensor_A X,Y) \\
  & \cong \RHom_A(X,\RHom_A(\RGammam A,Y)) \\
  & \stackrel{\rm (e)}{\cong} \RHom_A(X,\RHom_A(\RGammamopp A,Y)) \\
  & \stackrel{\rm (f)}{\cong} \RHom_A(X,Y).
\end{align*}
Here (d) is by \cite[thm.\ 1.6]{PJGorHom} again,
while (e) is by \cite[cor.\ 4.8]{VdBExist} and (f) is by the above
computation.
\end{proof}

\begin{Lemma}
\label{lem:standard_spectral_sequence}
For $X$ in $\D^-(\Gr\: A)$ and $Y$ in $\D^+(\Gr\: A)$ there is a convergent
spectral sequence
\[
  E_2^{mn} = \Ext_A^m(\h^{-n}\!X,Y) \Rightarrow \Ext_A^{m+n}(X,Y).
\]
\end{Lemma}

\begin{proof}
Let $J$ be an injective resolution of $Y$.  Consider the double
complex given by
\[
  M^{mn} = \Hom_A(X^{-m},J^n).
\]
The spectral sequence arising from the second standard filtration 
of the total complex $\Tot M$ gives the lemma's spectral sequence. 
\end{proof}

\begin{Lemma}
\label{lem:local_cohomology_spectral_sequence}
For $X$ in $\D^-(\Gr\: A)$ there is a convergent spectral sequence
\[
  E_2^{mn} = \Tor^A_{-m}(\HGammamopp^n\!A,X) \Rightarrow 
  \HGammam^{m+n}\!X.
\]
\end{Lemma}

\begin{proof}
Let $F$ be a flat resolution of $X$.  Consider the double complex
given by
\[
  M^{mn} = (\RGammamopp A)^m \otimes_A F^n.
\]
The spectral sequence arising from the second standard filtration
of $\Tot M$ gives the lemma's spectral sequence.

To see that the sequence has the indicated limit, one needs the
computation 
\[
  \Tot M \cong (\RGammamopp A) \LTensor_A X
  \stackrel{(a)}{\cong} (\RGammam A) \LTensor_A X
  \stackrel{(b)}{\cong} \RGammam X,
\]
where (a) is by \cite[cor.\ 4.8]{VdBExist} and (b) is by 
\cite[thm.\ 1.6]{PJGorHom}.
\end{proof}

\section{Two notions of regularity}
\label{sec:regularity}

The following is almost the classical definition of
Ca\-stel\-nu\-o\-vo-Mum\-ford regularity of graded modules, given over
polynomial algebras in \cite[dfn.\ p.\ 95]{EG} and more generally in
\cite[dfn.\ 2.1]{PJCM}.

\begin{Definition}[Castelnuovo-Mumford regularity]
\label{dfn:CMreg}
The complex $X$ in $\D(\Gr\: A)$ is called $p$-regular if
\[
  \HGammam^m(X)_{\geq p+1-m} = 0
\]
for all $m$.

If $X$ is $p$-regular but not $(p-1)$-regular, then I define the
{\em Ca\-stel\-nu\-o\-vo-Mum\-ford regularity of $X$} to be
\[
  \CMreg X = p.
\]
If $X$ is not $p$-regular for any $p$, then $\CMreg X = \infty$.
If $X$ is $p$-regular for every $p$ (that is, if $\HGammam(X) = 0$),
then $\CMreg X = -\infty$.
\end{Definition}

The following is the competing definition of Castelnuovo-Mumford
regularity given in \cite{AE}. In order not to confuse things, I have
to use a different name.

\begin{Definition}[Ext-regularity]
\label{dfn:Extreg}
The complex $X$ in $\D(\Gr\: A)$ is called $r$-Ext-regular if
\[
  \Ext_A^m(X,k)_{\leq -r-1-m} = 0
\]
for all $m$.

If $X$ is $r$-Ext-regular but not $(r-1)$-Ext-regular, then I define
the {\em Ext-regularity of $X$} to be
\[
  \Extreg X = r.
\]
If $X$ is not $r$-Ext-regular for any $r$, then $\Extreg X = \infty$. 
If $X$ is $r$-Ext-regular for every $r$ (that is, if $\Ext_A(X,k) = 0$),
then $\Extreg X = -\infty$.
\end{Definition}

\begin{Observation}
\label{obs:CMreg}
Let $X$ in $\D^{\b}_{\fg}(\Gr\: A)$ have $X \not\cong 0$.
Since $A$ has a balanced dualizing complex, the local duality theorem
\cite[thm.\ 4.18]{YDual} holds, so $\RGammam(X)^{\prime}$
is in $\D^{\b}_{\fg}(\Gr\: A^{\opp})$ and has $\RGammam(X)^{\prime}
\not\cong 0$. Hence $\CMreg X \not= \pm\infty$.

By \cite[cor.\ 4.8]{VdBExist} I have $\HGammam^n\!A \cong
\HGammamopp^n\!A$ for each $n$, whence
\[
  \CMreg({}_{A}A) = \CMreg(A_A).
\]
I denote this number by $\CMreg A$.
\end{Observation}

\begin{Observation}
\label{obs:Extreg}
Let $X$ in $\D^{\b}_{\fg}(\Gr\: A)$ have $X \not\cong 0$. It is easy
to see $\Ext_A(X,k) \not\cong 0$ whence $\Extreg X \not=
-\infty$. However, $\Extreg X = \infty$ is possible.

If $F$ is a minimal free resolution of $X$, then $X$ is
$r$-Ext-regular exactly if the generators of $F_m$ are placed in
degrees less than or equal to $r+m$ for each $m$.  This has a nice
consequence: From considering $\Tor^A(k_A,{}_{A}k)$ it follows that
the minimal free resolutions of $k_A$ and ${}_{A}k$ have their
generators placed in the same degrees. Hence
\[
  \Extreg(k_A) = \Extreg({}_{A}k).
\]
I denote this number by $\Extreg k$.
\end{Observation}

%\begin{Remark}
%Castelnuovo-Mumford regularity has the advantage over Ext-regularity
%that if $X$ is in $\D_{\fg}^{\b}(\Gr\: A)$, then $\CMreg X$ is a
%finite number, whereas $\Extreg X$ might well be infinite.
%
%On the other hand, Ext-regularity has the advantage over
%Castelnuovo-Mumford regularity that it encodes a more tangible piece
%of information, namely the behaviour of a minimal free resolution.
%\end{Remark}

The following two theorems show that the notions of
Castelnuovo-Mumford and Ext-regularity enjoy a close
relationship. Note the structural similarity between the proofs.

\begin{Theorem}
\label{thm:main1}
Given $X$ in $\D^{\b}_{\fg}(\Gr\: A)$ with $X \not\cong 0$. Then
\[
  \Extreg X  \leq  \CMreg X + \Extreg k.
\]
\end{Theorem}

\begin{proof}
Observation \ref{obs:CMreg} gives $\CMreg X \not= -\infty$, so for
$\Extreg k = \infty$ the theorem makes sense and holds trivially. So I
can assume that $\Extreg k = r$ is finite. By observation
\ref{obs:Extreg}, the minimal free resolution $F$ of $k_A$ then has
the generators of $F_m$ placed in degrees less than or equal to $r+m$
for each $m$, so $F_m$ can be written as a finite coproduct
\[
  F_m = \coprod_j A(-\sigma_{mj})
\]
with $\sigma_{mj} \leq r+m$. Taking Matlis duals, $I = F^{\prime}$ is a
minimal injective resolution of ${}_{A}k$ which has
\[
  I^m = \coprod_j A^{\prime}(\sigma_{mj}),
\]
still with 
\begin{equation}
\label{equ:sigma}
  \sigma_{mj} \leq r+m.
\end{equation}

Set $p = \CMreg X$ and $Z = \RGammam X$. Then
\[
  \h^{-n}(Z)_{\geq p+1+n} = \h^{-n}(\RGammam X)_{\geq p+1+n} =
  \HGammam^{-n}(X)_{\geq p+1+n} = 0
\]
for each $n$ whence
\begin{equation}
\label{equ:Zprime}
  ((\h^{-n}\!Z)^{\prime})_{\leq -p-1-n} = 0.
\end{equation}
Now, $\Ext_A^m(\h^{-n}\!Z,k)$ is a subquotient of
\begin{align*}
  \Hom_A(\h^{-n}\!Z,I^m)
  & = \Hom_A(\h^{-n}\!Z,\coprod_j A^{\prime}(\sigma_{mj})) \\
  & \cong \coprod_j (\h^{-n}\!Z)^{\prime}(\sigma_{mj}),
\end{align*}
and this vanishes in degrees less than or equal to $-p-1-n-r-m$ by
equations \eqref{equ:sigma} and \eqref{equ:Zprime}, so also
\begin{equation}
\label{equ:Ext}
  \Ext_A^m(\h^{-n}\!Z,k)_{\leq -p-1-n-r-m} = 0.
\end{equation}

Lemma \ref{lem:standard_spectral_sequence} provides a convergent
spectral sequence
\[
  E_2^{mn} = \Ext_A^m(\h^{-n}\!Z,k) \Rightarrow \Ext_A^{m+n}(Z,k),
\]
and since equation \eqref{equ:Ext} shows $(E_2^{mn})_{\leq -p-1-r-(m+n)} =
0$ it follows that
\begin{equation}
\label{equ:Ext2}
  \Ext_A^q(Z,k)_{\leq -p-1-r-q} = 0
\end{equation}
for each $q$.

Finally, proposition \ref{prp:RHom_isomorphism} gives 
\[
  \Ext_A^q(Z,k) = \Ext_A^q(\RGammam X,k) \cong \Ext_A^q(X,k),
\]
so equation \eqref{equ:Ext2} implies
\[
  \Ext_A^q(X,k)_{\leq -p-1-r-q} = 0
\]
for each $q$, showing
\[
  \Extreg X  \leq  p+r  =  \CMreg X + \Extreg k.
\]
\end{proof}

\begin{Theorem}
\label{thm:main2}
Given $X$ in $\D^{\b}_{\fg}(\Gr\: A)$ with $X \not\cong 0$. Then
\[
  \CMreg X  \leq  \Extreg X + \CMreg A.
\]
\end{Theorem}

\begin{proof}
Observation \ref{obs:CMreg} gives $\CMreg A \not= -\infty$, so for
$\Extreg X = \infty$ the theorem makes sense and holds trivially. So I
can assume that $\Extreg X = r$ is finite. By observation
\ref{obs:Extreg}, the minimal free resolution $F$ of $X$ then has
the generators of $F_m$ placed in degrees less than or equal to $r+m$
for each $m$, so $F_m$ can be written as a finite coproduct
\[
  F_m = \coprod_j A(-\sigma_{mj})
\]
with
\begin{equation}
\label{equ:sigma2}
  \sigma_{mj} \leq r+m.
\end{equation}

Set $p = \CMreg A$. Observation \ref{obs:CMreg} gives $\CMreg(A_A)
= \CMreg A$, so I get
\begin{equation}
\label{equ:HGammamopp}
  \HGammamopp^n(A)_{\geq p+1-n} = 0
\end{equation}
for each $n$.

Now, $\Tor^A_{-m}(\HGammamopp^n\!A,X)$ is a subquotient of
\begin{align*}
  \HGammamopp^n(A) \otimes_A F_{-m}
  & \cong \HGammamopp^n(A) \otimes_A \coprod_j A(-\sigma_{-m,j}) \\
  & \cong \coprod_j \HGammamopp^n(A)(-\sigma_{-m,j}),
\end{align*} 
and this vanishes in degrees larger than or equal to $p+1-n+r-m$ by
equations \eqref{equ:sigma2} and \eqref{equ:HGammamopp}, so also
\begin{equation}
\label{equ:Tor}
  \Tor^A_{-m}(\HGammamopp^n\!A,X)_{\geq p+1-n+r-m} = 0.
\end{equation}

Lemma \ref{lem:local_cohomology_spectral_sequence} provides a convergent
spectral sequence
\[
  E_2^{mn} = \Tor^A_{-m}(\HGammamopp^n\!A,X) \Rightarrow 
  \HGammam^{m+n}\!X,
\]
and since equation \eqref{equ:Tor} shows $(E_2^{mn})_{\geq p+1+r-(m+n)} =
0$, it follows that
\begin{equation}
\label{equ:Tor2}
  \HGammam^q(X)_{\geq p+1+r-q} = 0
\end{equation}
for each $q$, showing
\[
  \CMreg X  \leq  p+r  =  \Extreg X + \CMreg A.
\]
\end{proof}

Let me end the section with some easy consequences.  First recall the
following definition.

\begin{Definition}
\label{def:Koszul}
The algebra $A$ is called {\em Koszul} if $\Extreg k = 0$.
\end{Definition}

For $A$ to be Koszul means exactly that the minimal free resolutions
of ${}_{A}k$ and $k_A$ are linear, cf.\ observation \ref{obs:Extreg}. 

The following corollary is immediate from theorems \ref{thm:main1} and
\ref{thm:main2}.

\begin{Corollary}
\label{cor:regs_coincide}
Suppose that $A$ is Koszul and has $\CMreg A = 0$. Then any $X$ in
$\D^{\b}_{\fg}(\Gr\: A)$ has $\Extreg X = \CMreg X$.
\end{Corollary}

The following corollary is also immediate from observation
\ref{obs:CMreg} and theorem \ref{thm:main1}.  It extends
\cite[thm.\ 1]{AE} and \cite[thm.\ 1]{AP} to the
non-com\-mu\-ta\-ti\-ve case. 

\begin{Corollary}
\label{cor:Extreg_finite}
Suppose that $A$ has $\Extreg k < \infty$. Then any $X$ in
$\D^{\b}_{\fg}(\Gr\: A)$ has $\Extreg X < \infty$.
\end{Corollary}

\section{Linear free resolutions}
\label{sec:linear}

The following main result is a simultaneous extension of
\cite[cor.\ 2]{AE} (to the non-commutative case) and \cite[thm.\
2.6]{PJCM} (to the non-AS-regular case). 

Recall that $A$ is the algebra of setup \ref{set:blanket}.

\begin{Theorem}
\label{thm:linear}
Suppose that $A$ is Koszul, and let $M$ in $\Gr\: A$ be finitely
generated with $M \not\cong 0$.  Then for $s \geq \CMreg M$, the
minimal free resolution of $M_{\geq s}(s)$ is linear.  (Note that
$\CMreg M$ is finite.)
\end{Theorem}

\begin{proof} 
The result clearly holds if $M_{\geq s}(s)$ is $0$, so I can assume
$M_{\geq s}(s) \not\cong 0$.

Let $F$ be the minimal free resolution of $M_{\geq s}(s)$. As $M_{\geq
s}(s)$ sits in non-negative degrees, it is clear for each $m$ that
$F_m$ has no generators placed in degrees strictly smaller than
$m$. Hence it is enough to prove for each $m$ that $F_m$ also has no
generators placed in degrees strictly larger than $m$. By observation
\ref{obs:Extreg} this is the same as proving
\begin{equation}
\label{equ:Extreg_condition}
  \Extreg(M_{\geq s}(s)) \leq 0. 
\end{equation}

Since $A$ is Koszul, $\Extreg k = 0$ holds. By theorem
\ref{thm:main1}, the inequality \eqref{equ:Extreg_condition} will
therefore follow from $\CMreg(M_{\geq s}(s)) \leq 0$, which is again the
same as $\CMreg(M_{\geq s}) \leq s$, that is
\[
  \HGammam^m(M_{\geq s})_{\geq s+1-m} = 0
\]
for each $m$. To show this is easy:

There is a short exact sequence $0 \rightarrow M_{\geq s}
\longrightarrow M \longrightarrow M/M_{\geq s} \rightarrow 0$
resulting in a long exact sequence consisting of pieces
\[
  \HGammam^m(M_{\geq s}) \longrightarrow
  \HGammam^m(M) \longrightarrow 
  \HGammam^m(M/M_{\geq s}).
\] 
Now combine this with $\HGammam^m(M)_{\geq s+1-m} = 0$ for each $m$
(because I have $s \geq \CMreg M$) and
\[
  \HGammam^m(M/M_{\geq s}) \cong
  \left\{
    \begin{array}{cl}
      M/M_{\geq s} & \mbox{for $m = 0$,} \\
      0            & \mbox{for $m \geq 1$}
    \end{array}
  \right.
\]
(because $M/M_{\geq s}$ is torsion).
\end{proof}

\end{document}